\documentclass[12pt]{amsart}
\usepackage{graphicx}
\usepackage{amssymb}
\usepackage{amsfonts}
\usepackage{amsmath}
\usepackage{array}
\usepackage{rotating}
\usepackage[matrix,frame]{xypic}

\def\FQSym{{\bf FQSym}}

\def\std{{\rm std}}

\def\Eul{{\mathcal E}}
\def\<{\langle}
\def\>{\rangle}

\def\F{{\bf F}}

\def\G{{\bf G}}

\def\SG{{\mathfrak S}}

\def\A{{\bf A}}

\def\Sym{{\bf Sym}}

\def\Tabvrule{\vrule width-0.4pt}       
\def\Tabhrule{\hrule \hrule height-0.4pt} 
\def\Tabstrut{\vrule height2.2ex 
                     depth0.8ex  
                     width0ex    
\relax}

\def\PasCase#1{\omit%
            $\vcenter{\hbox {\vbox to 0.4pt{}}
               \hbox{\makebox[3ex]{\Tabstrut$#1$}}}%
               \Tabvrule$}
\def\PasCasePoint{\PasCase{\cdot}}
\def\DessinCarre#1{%
    \vcenter{\hbox{}\hrule
             \hbox{\vrule\makebox[3ex]{\Tabstrut$#1$}\vrule}\Tabhrule}%
             \Tabvrule}
\def\GenRuban#1{\vcenter{\halign{&$\DessinCarre{##}$\cr#1}}\egroup}

\def\sTabvrule{\vrule width-0.4pt}
\def\sTabhrule{\hrule \hrule height-0.4pt}
\def\sTabstrut{\vrule height1.6ex depth0.6ex width0ex \relax}
\def\sDessinCarre#1{%
    \vcenter{\hbox{}\hrule
             \hbox{\vrule\makebox[2.3ex]%
                  {\sTabstrut$\scriptstyle#1$}\vrule}\sTabhrule}%
             \sTabvrule}
\def\sGenRuban#1{\vcenter{\halign{&$\sDessinCarre{##}$\cr#1}}\egroup}

\def\ruban{%
  \bgroup
  \let\ =\omit
  \let\\=\cr
  \let\x=\times
  \let\.=\PasCasePoint
  \offinterlineskip
  \GenRuban}

\def\sruban{%
  \bgroup
  \let\ =\omit
  \let\x=\times
  \let\\=\cr
  \offinterlineskip
  \sGenRuban}

\title{Noncommutative Symmetric Functions and an Amazing Matrix}

\author[ J.-C.~Novelli and J.-Y.~Thibon]%
{Jean-Christophe Novelli and Jean-Yves Thibon}

\address[] {Laboratoire d'Informatique  Gaspard Monge, Universit\'e Paris-Est Marne-la-Vall\'ee \\
5 Boulevard Descartes \\Champs-sur-Marne \\77454 Marne-la-Vall\'ee cedex 2 \\
FRANCE}
\email[Jean-Christophe Novelli]{novelli@univ-mlv.fr}
\email[Jean-Yves Thibon]{jyt@univ-mlv.fr} 
\date{}

\begin{document}

\begin{abstract}
We present a simple way to derive the results of Diaconis and Fulman [arXiv:1102.5159]
in terms of noncommutative symmetric functions.
\end{abstract}

\maketitle


\section{Introduction}

In \cite{DF1,DF2}, Diaconis and Fulman investigate a remarkable family of matrices
introduced by Holte \cite{Hol} in his analysis of the process of ``carries'' in the
addition of random integers in base $b$. 

The aim of this note is to show that the results of \cite{DF2} can be derived
in a simple and natural way within the formalism of noncommutative symmetric
functions \cite{NCSF1}.

This is possible thanks to the following equivalent characterization of
the ``amazing matrix'' $P$ (Theorem 2.1 of \cite{DF1}):
\begin{quote}
The number of descents in successive $b$-shuffles of $n$ cards form a Markov
chain on $\{0,1,\ldots,n-1\}$ with transition matrix $P(i,j)$.
\end{quote}

Such random processes involving descents of permutations can usually be interpreted
in the descent algebra of the symmetric group. Here, it is only the number of descents
which is involved, so that one can in fact compute in the Eulerian subalgebra.

We assume that the reader is familiar with the notations of \cite{NCSF1}.

\section{The Eulerian algebra}

This is a commutative subalgebra of dimension $n$ of the group algebra of
the symmetric group $\SG_n$, and in fact of its descent algebra $\Sigma_n$. It was apparently
first  introduced in \cite{BMP} under the name {\em algebra of permutors}\footnote{
A self-contained and elementary presentation of the main results of \cite{BMP}
can be found in \cite{NCSF1}.}.
It is spanned by the Eulerian idempotents, or, as well, by the sums of
permutations having the same number of descents. 

It is easier to work with all symmetric groups at the same time, with the
help of generating functions. Recall that the algebra of noncommutative
symmetric functions $\Sym$ is endowed with an internal product $*$,
for which each homogeneous component $\Sym_n$ is anti-isomorphic to $\Sigma_n$.

Recall also the following definitions from \cite{NCSF1}.
We denote by $\sigma_t$ or $\sigma_t(A)$ the generating series of the
complete symmetric functions $S_n$:
\begin{equation}
\sigma_t(A)=\sum_{n\ge 0}t^n S_n(A).
\end{equation}
The {\em Eulerian idempotents} $E_n^{[k]}$ are the homogenous components
of degree $n$ in the series $E^{[k]}$ defined by
\begin{equation}
\sigma_t(A)^x=\sum_{k\ge 0}x^k E^{[k]}(A).
\end{equation}
We have
\begin{equation}
E_n^{[k]}*E_n^{[l]} = \delta_{kl}E_n^{[k]}\,,\quad\text{and}\quad\sum_{k=1}^nE_n^{[k]}=S_n,
\end{equation}
so that the $E_n^{[k]}$ span a commutative $n$-dimensional $*$-subalgebra of
$\Sym_n$, denoted by $\Eul_n$ and called the Eulerian subalgebra.

The {\em noncommutative Eulerian polynomials} are defined by
\begin{equation}
{\mathcal A}_n(t) =
\ \sum_{k=1}^n \ t^k\, \Big(
\sum_{{\scriptstyle |I|=n}\atop{\scriptstyle \ell(I)=k}} R_I \,
\Big)
=
\ \sum_{k=1}^n \ {\bf A}(n,k)\, t^k \, ,
\end{equation}
where $R_I$ is the ribbon basis.
The following facts can be found (up to a few misprints) in \cite{NCSF1}.
The generating series of the ${\mathcal A}_n(t)$ is
\begin{equation}
{\mathcal A}(t) := \ \sum_{n\ge 0} \, {\mathcal A}_n(t)
=
(1-t) \, \left( 1 - t\, \sigma_{1-t} \right)^{-1} \ .
\end{equation}

\medskip
Let ${\mathcal A}_n^*(t) = (1-t)^{-n}\, {\mathcal A}_n(t)$.
Then,
\begin{equation}
{\mathcal A}^*(t)
:=
\ \sum_{n\ge 0} \, {\mathcal A}_n^*(t)
=
\sum_{I} \
\left( \displaystyle {t \over 1-t} \right)^{\ell(I)} \, S^I \ .
\end{equation}
This last formula can also be written in the form
\begin{equation} \label{GEN*}
{\mathcal A}^*(t)
=
\ \sum_{k\ge 0} \ \left(
{t\over 1-t}\right)^k \left( S_1+S_2+S_3+\cdots\, \right)^k
\end{equation}
or
\begin{equation}\label{GEN_A}
{1\over 1-t\, \sigma_1(A)}
=
\ \sum_{n\ge 0}\ {{\mathcal A}_n(t)\over (1-t)^{n+1}} \ .
\end{equation}
Let $S^{[k]}=\sigma_1(A)^k$ be the coefficient of $t^k$ in this series. In degree $n$,
\begin{equation}\label{S2E}
S_n^{[k]}=\sum_{I\vDash n, \ell(I)\le k}{k\choose \ell(I)}S^I=\sum_{i=1}^nk^iE_n^{[i]}\,.
\end{equation}
This is another basis of $\Eul_n$. 
Expanding the factors $(1-t)^{-(n+1)}$ in the right-hand side
of (\ref{GEN_A}) by the binomial
theorem, and taking the coefficient of $t^k$ in the term of
weight $n$ in both sides, we get
\begin{equation}
S_n^{[k]}
=
\ \sum_{i=0}^k \ {n+i \choose i} \, {\bf A}(n,k-i) \, .
\end{equation}
Conversely,
\begin{equation}
{ {\mathcal A}_n(t) \over (1-t)^{n+1} }
=
\ \sum_{k\ge 0}\ t^k\, S_n^{[k]} \ ,
\end{equation}
so that
\begin{equation}\label{A2S}
{\bf A}(n,p)
=
\ \sum_{i=0}^p\ (-1)^i\, {n+1\choose i}\, S_n^{[p-i]} \ .
\end{equation}
%
%
The expansion of the $E_n^{[k]}$
on the basis ${\bf A}(n,i)$, which is a noncommutative analog
of Worpitzky's identity (see \cite{Ga} or \cite{Lod89}) is
\begin{equation}\label{Wor}
\sum_{k=1}^n \ x^k\, E_n^{[k]}
=
\ \sum_{i=1}^n \ {x+n-i\choose n}\, {\bf A}(n,i) \ .
\end{equation}
Indeed, when $x$ is a positive integer $N$,
\begin{equation}
\sum_{k=1}^n \ N^k\, E_n^{[k]}
= S_n(NA) = \sum_{I\vDash n}F_I(N)R_I(A)
\end{equation} 
where $F_I$ are the fundamental quasi-symmetric functions,
and for a composition $I=(i_1,\ldots,i_r)$ of $n$,
\begin{equation}
F_I(N)={N+n-r\choose n}\,.
\end{equation}


\section{The $b$-shuffle process}

For a positive integer $b$, the $b$-shuffle permutations in $\SG_n$ are
the inverses of the permutations with at most $b-1$ descents. Thus, the
$b$-shuffle operator can be identified with 
$S_n^{[b]}$ ({\it i.e.}, with $*$-multiplication by
$S_n^{[b]}$). It belongs to the Eulerian
algebra, so that it preserves it, and it makes sense to compute its matrix
in the basis $\A(n,k)$. Note that since $\Eul_n$ is commutative, it does not
matter whether we multiply on the right or on the left.

Note that the $b$-shuffle process is an example of what Stanley has
called the QS-distribution \cite{St}. It is the probability distribution
on permutations derived by assigning probability $b^{-1}$ to the
first $b$ positive integers, see \cite{NCSF6} for a simplified version.

Summarizing, we want to compute the coefficients $P_{ij}(b)$ defined by
\begin{equation}
S_n^{[b]} * \A(n,j) = \sum_{i=1}^n P_{ij}(b) \A(n,i)\,.
\end{equation}
From (\ref{S2E}), it is clear that
\begin{equation}
S_n^{[p]}*S_n^{[q]}=S_n^{[pq]}
\end{equation}
so that, using (\ref{A2S}), we obtain
\begin{equation}
\begin{split}
S_n^{[b]} * \A(n,j) &= \sum_{r=0}^j(-1)^r{n+1\choose r}S_n{[b(j-r)]}\\
&= \sum_{r=0}^j(-1)^r{n+1\choose r}\sum_{k=0}^{b(j-r)}{n+k\choose k}\A(n,b(j-r)-k)\,.
\end{split}
\end{equation}
The coefficient of $\A(n,i)$ in this expression is therefore
\begin{equation}
P_{ij}(b)=\sum_{r=0}^j(-1)^r{n+1\choose r}{n+b(j-r)-i \choose n}
\end{equation}
These are the coefficients of the amazing matrix (up to a shift of 1 on the indices
$i,j$, and a global normalization factor $b^n$ so as the probabilities sum up to 1).

Since the $E_n^{[k]}$ form a basis of orthogonal idempotents in $\Eul_n$,
it is reasonable to introduce a scalar product such that
\begin{equation}
\<E_n^{[i]}|E_n^{[j]}\> = \delta_{ij}\,.
\end{equation}
Then, the $b$-shuffle operator is self-adjoint. 
Its orthonormal basis of eigenvectors is clearly $E_n^{[k]}$
(with eigenvalues $b^k$).

In terms of coordinates,
since we are working in the
non-orthogonal basis $\A(n,i)$,
its right eigenvector of eigenvalue $b^j$ is the column vector
whose $i$th component is the coefficient of $E_n^{[j]}$ on $\A(n,i)$,
that is, the coefficient of $x^j$ in ${x+n-i\choose n}$, thanks to
(\ref{Wor}). By duality, its left eigenvector associated with the 
eigenvalue $b^i$ is the row vector whose $j$th component is 
\begin{equation}
\<\A(n,j)|E_n^{[i]}\> = \sum_{r=0}^j(-1)^r{n+1\choose r}(j-r)^i.
\end{equation}
This is precisely the Foulkes character table (up to indexation,
the Frobenius characteristic of $\chi^{n,k}$ is the commutative image\footnote{
These commutative symmetric functions have been studied in \cite{Desar}.} of $\A(n,n-k)$).

\section{Other examples}

\subsection{Determinant of the Foulkes character table}

This is the determinant of the matrix $F$
\begin{equation}
F(i,j) = \<\A(n,i), E_n^{[j]}\>\,\quad i,j=1,\ldots,n\,.
\end{equation}
Because of the triangularity property
\begin{equation}
\A(n,i) = S_n^{[i]} + \sum_{r=1}^i(-1)^r{n+1\choose r}S_n^{[i-r]},
\end{equation}
we have as well
\begin{equation}
\det F = \det G\,\quad\text{where}\ G(i,j)=\<S_n^{[i]},E_n^{[j]}\> = i^j
\end{equation}
a Vandermonde determinant which evaluates to $n!(n-1)!\cdots 2!1!$.

\subsection{Descents of $b^r$-riffle shuffles}

Recall from \cite{NCSF6} that $\FQSym$ is
an algebra based on all permutations and that it has two bases
\begin{equation}
\G_\sigma= \sum_{\std(w)=\sigma} w = \F_{\sigma^{-1}}
\end{equation}
which are mutually adjoint for its natural scalar product
\begin{equation}
\<\F_\sigma, \G_\tau\>=\delta_{\sigma,\tau}\,.
\end{equation}
Under the embedding of $\Sym$ into 
$\FQSym$, the $b^r$-shuffle operator is
\begin{equation}
(S_n^{[b]})^{*r} = S_n^{[b^r]} = \sum_{\sigma\ b^r-\text{shuffle}}\F_\sigma\,.
\end{equation}
The generating function of $b^r$-shuffle by number of descents is therefore
its scalar product in $\FQSym$ with the noncommutative Eulerian polynomial
\begin{equation}
{\mathcal A}_n(t)=\sum_{k=1}^nt^k \A(n,k) = \sum_{\tau\in\SG_n}t^{d(\tau)+1}\G_\tau .
\end{equation}
Recall that
\begin{equation}
{\mathcal A}_n(t)=(1-t)^{n+1}\sum_{k=1}^nt^k S_n^{[k]}
\end{equation}
so that
\begin{equation}
\<S_n^{[b^r]},{\mathcal A}_n(t)\>=(1-t)^{n+1}\sum_{k=1}^nt^k \<S_n^{[b^r]},S_n^{[k]}\>
\end{equation}
Now, when one factor $P$ of a scalar product $\<P,Q\>$ in $\FQSym$ is in $\Sym$,
one has $\<P,Q\>=\<p,Q\>$ where $p=\underline{P}$ is the commutative
image of $P$ in $QSym$, and the bracket is now the duality between $\Sym$
and $QSym$. Furthermore, when $p$ in in $Sym$, then, the scalar product reduces
to $\<p,q\>$, where $q=\underline{Q}$ is the commutative image of $Q$
in $Sym$, and the bracket is now the ordinary scalar product of symmetric
functions (see \cite{NCSF6}).
Thus,
\begin{equation}
 \<S_n^{[b^r]},S_n^{[k]}\> = \<h_n(b^rX),h_n(kX)\>=h_n(b^rk)={b^rk+n-1\choose n}
\end{equation}
($\lambda$-ring notation) and we are done:
\begin{equation}
\<S_n^{[b^r]},{\mathcal A}_n(t)\>=(1-t)^{n+1}\sum_{k=1}^nt^k {b^rk+n-1\choose n}\,.
\end{equation}


\end{document}